\documentclass[12pt]{amsart}
\usepackage{amsmath, amsthm,url, amssymb,setspace,epsfig, hyperref, mathrsfs}
\usepackage[all]{xy}
\SelectTips{cm}{}

\newtheorem{thm}{Theorem}[section]
\newtheorem{prop}[thm]{Proposition}
\newtheorem{cor}[thm]{Corollary}
\newtheorem{lemma}[thm]{Lemma}

\theoremstyle{definition}
\newtheorem{defe}[thm]{Definition}
\theoremstyle{remark}
\newtheorem{exam}[thm]{Example}

\newtheorem{ntn}[thm]{Convention}

\newtheorem{rem}[thm]{Remark}

\numberwithin{equation}{section}

\newcommand         {\rar}[1]       {\stackrel{#1}{\longrightarrow}}
\newcommand         {\isom}         {\rar{\simeq}}
\newcommand         {\eqv}          {\rar{\sim}}
\newcommand         {\commentout}[1]    {}

\def\inv{^{-1}}

\def\qm{\mathrm{qm}}

\def    \ra     {\rightarrow}
\def    \inj    {\hookrightarrow}
\def    \onto   {\twoheadrightarrow}
\def    \*      {\times}
\def    \ten    {\otimes}

\def    \SP     {\mathrm{Sp}}

\def    \Exp    {\mathrm{Exp}}
\def    \Log    {\mathrm{Log}}
\def    \CH     {\mathrm{CH}}
\def    \Nilp   {\mathrm{Nilp}}
\def    \Hom    {\mathrm{Hom}}
\def    \Ker    {\mathrm{Ker}}
\def    \ind    {\mathrm{ind}}
\def    \Ad     {\mathrm{Ad}}

\def    \th     {\mathrm{th}}
\def    \Z      {\mathrm{Z}}
\def    \Id     {\mathrm{Id}}
\def    \ph     {\varphi}

\def    \Hom    {\mathrm{Hom}}
\def    \BQ     {\mathrm{BQ}}

\def    \bZm    {\bZ/{m\bZ}}

\def    \Q      {\mathrm{Q}}
\def    \B      {\mathrm{B}}
\def    \K      {\mathrm{K}}
\def    \det    {\mathrm{det}}
\def    \Det    {\mathrm{Det}}
\def    \is     {\mathrm{is}}
\def    \th    {\mathrm{th}}

\def    \QM     { \Q \cM }
\def    \BQM    {\B \Q \cM}
\def    \PKM    {\V( \cM)}
\def    \V      {\mathrm{V}}
\def    \VM     {\V (\cM)}

\def	\gr		{\mathrm{gr}}

\def    \fz     {\mathfrak{z}}

\def    \fg     {\mathfrak{g}}

\def    \fp     {\mathfrak{p}}
\def    \fS     {\mathfrak{S}}
\def    \fR     {\mathfrak{R}}
\def    \fq     {\mathfrak{q}}
\def    \fh     {\mathfrak{h}}

\def    \fl     {\mathfrak{l}}
\def    \nilp  {\mathfrak{nilp}}
\def    \fD     {\mathfrak{D}}
\def    \DPic   {\fD_{\Pic}}

\def    \bZ     {\mathbb{Z}}

\def    \bC     {\mathbb{C}}
\def    \bCt    {\bC^\times}
\def    \Fq     {\mathbb{F}_q}
\def    \Fp     {\mathbb{F}_p}
\def    \Fpt    {\Fp^\times}
\def    \Fptt   {{\Fpt}/{(\Fpt)^2}}

\def    \cE     {\mathcal{E}}
\def    \cP     {\mathcal{P}}
\def    \cMo    {({\cM,\is})}
\def    \cS     {\mathcal{S}}
\def    \cM     {\mathcal{M}}

\def    \bf     {\bar{f}}

\def    \bA     {\bar{A}}
\def    \bfp    {{\fp}^f}
\def    \bfg    {{\fg}^f}
\def    \bo     {\bar{o}}

\def    \Vect       {\mathbf{Vect}}
\def    \Pic        {\mathbf{Pic}}

\def    \Abp     {\mathbf{Ab}_p}
\def    \Vectp    {\mathbf{Vect}_{\Fp}}

\def    \bW     {\mathbb{W}}
\def    \Fp     {\mathbb{F}_p}
\def    \bZ     {\mathbb{Z}}

\def    \tL     {\tilde{L}}

\def    \tfp    {\tilde{\mathfrak{p}}}
\def    \tfg    {\tilde{\fg}}

\DeclareMathOperator\len{len}

\begin{document}
\title[Compatible Intertwiners for Finite Nilpotent Groups]{Compatible Intertwiners for Representations of Finite Nilpotent Groups}

\author{Masoud Kamgarpour}
\address{The University of British Columbia,  Vancouver, Canada V6T 1Z2}
\email{masoud@math.ubc.ca}

\author{Teruji Thomas}
\address{The University of Edinburgh,  Edinburgh, United Kingdom EH9 3JZ}
\email{t.thomas@ed.ac.uk}

\thanks{M.K. was supported by NSERC PDF grant. T.T. was supported by a JRF at Merton College, Oxford and a Seggie Brown Fellowship at Edinburgh. }


\bibliographystyle{alpha}

\begin{abstract} We sharpen the orbit method for finite groups of small nilpotence class by associating  representations to functionals on the corresponding Lie rings.  This amounts to describing compatible intertwiners between representations parameterized by an additional choice of polarization.  Our construction is motivated by the theory of the linearized Weil representation of the symplectic group. In particular, we provide generalizations of the Maslov index and the determinant functor to the  context of finite abelian groups.
\end{abstract}
\maketitle

{\small \setcounter{tocdepth}{1} \tableofcontents}

\section{Introduction}\label{s:introduction}

This paper is about the complex irreducible representations of finite nilpotent groups $G$ of small nilpotence class. Here `small' means `smaller than any prime divisor of $|G|$.' However, it is convenient to note that such a group canonically decomposes as the direct product of its Sylow subgroups, so that the whole theory reduces to the case where $G$ is a finite $p$-group of nilpotence class less than $p$. We will also assume that $p$ is an odd prime, since $p=2$ forces $G$ to be abelian.

\subsection{Statement of the problem} \label{ss:statement}
Let, then,   $G$ be a finite $p$-group of nilpotence class less than $p$. In \cite[\S2]{DB}, M. Boyarchenko and V. Drinfeld have modified Kirillov's orbit method ~\cite{Kirillov} to describe irreducible representations of $G$ (see also \cite[\S6]{Howe-Nilpotent}).  The construction works as follows (with more details in \S\S\ref{s:polarizations}--\ref{s:repOfpGroups}).
To $G$ one associates a {\it Lie ring}, that is to say, an abelian group $\fg$ equipped with a bracket $[\cdot\,,\cdot]$ satisfying the Jacobi identity and $[x,x]=0$ for all $x\in\fg$.  The underlying set of $\fg$ equals that of $G$.
The relationship between $\fg$ and $G$ is analogous to that of Lie algebras and Lie groups.
The multiplication in $G$ is given by the usual Campbell-Hausdorff formula:
\begin{equation}\label{eq:CampbellHaus}
x*y=x+y+\frac{[x,y]}{2}+(\mbox{higher order terms}).
\end{equation}
There is a one-to-one correspondence between subgroups of $G$ and subrings of $\fg$. Furthermore, one can define the adjoint action of $G$ on $\fg$ in an obvious way.

Pick $f$ in $\fg^*:=\Hom(\fg,\bCt)$. A {\it polarization} for $f$ is a subring $\fp\subseteq \fg$ satisfying $f([\fp,\fp])=\{1\}$, such that $\fp$  is maximal among all additive {\it subgroups} of $\fg$ with this property. Let $P$ be the subgroup of $G$ corresponding to a polarization $\fp$. It follows from (\ref{eq:CampbellHaus}) that $f$ defines a homomorphism $\chi_f\colon P\ra \bCt$. The orbit method states that $\rho_{f,\fp}:=\ind_P^G \chi_f$ is an irreducible representation of $G$. Furthermore, the isomorphism class of $\rho_{f,\fp}$ depends only on the co-adjoint orbit of $f$, and every isomorphism class arises in this way.

\subsubsection{}\label{sss:goal}
The main goal of this article is to sharpen the orbit method by constructing a canonical representation $\rho_f$ of $G$ for each $f\in\fg^*$, independent of the choice of polarization. We do this by constructing compatible intertwiners between the various representations $\rho_{f,\fp}$. The intertwiner between $\rho_{f,\fp_1}$ and $\rho_{f,\fp_2}$ depends on a choice of \emph{orientations} of $\fp_1$ and $\fp_2$. An orientation
 of a vector space $V$ is a non-zero element\footnote{In the literature it is usual, and in practice it suffices, to identify two orientations $o_1$ and $o_2$ if $o_1=a^2 o_2$ for some $a$ in the ground field.  In that sense, for example,
any real vector space has two orientations. In \cite{RonnyCanonical}, the less picturesque term `enhancement' is used in place of what we call orientation. } of the determinant (highest exterior power) $\det(V)$; in Appendix \ref{s:Det}, we develop the corresponding theory of determinants of finite abelian $p$-groups.  We can then define an oriented polarization to be a pair $\tfp=(\fp,o)$ consisting of a polarization $\fp$ and a nonzero element $o\in \det(\fp)$.

Thus, for every pair $\tfp_1,\tfp_2$ of oriented polarizations, we describe intertwining operators $\Psi_{\tfp_1,\tfp_2}\in \Hom_G(\rho_{f,\fp_2},\rho_{f,\fp_1})$ that are \emph{compatible}; i.e., for any three oriented polarizations $\tfp_1,\tfp_2,\tfp_3$,  the following diagram commutes
\begin{equation}
\xymatrix{ \rho_{f,\fp_3} \ar[rr]^{\Psi_{\tfp_1,\tfp_3}} \ar[dr]_{\Psi_{\tfp_2,\tfp_3}} & & \rho_{f,\fp_1}. \\
&   \rho_{f,\fp_2} \ar[ur]_{\Psi_{\tfp_1,\tfp_2}}& }
\end{equation}
Let $\Lambda_f$ be the set of all oriented polarizations of $f$, and write $V_{f,\fp}$ for the representation space of $\rho_{f,\fp}$. Then the desired $\rho_f$ is the diagonal representation of $G$ on
$$V_f:=\{(v_{\tfp}) \in \bigoplus_{\tfp\in\Lambda_f} V_{f,\fp}\mid v_{\tfp_1}=\Psi_{\tfp_1,\tfp_2}(v_{\tfp_2})\mbox{ for all $\tfp_1,\tfp_2\in\Lambda_f$}\}.$$
The projection of $V_f$ to each $V_{f,\fp}$ is an isomorphism between representations of $G$.

\begin{rem} \label{r:HeisenbergGroupWeilRep} \emph{The Weil Representation.}
 Let $(A,\omega)$ be a symplectic vector space over $\Fq$. The Heisenberg group $H$ is the set $A\times \Fq$ equipped with the multiplication
\[
(v,t)*(v',t')=(v+v',t+t'+\frac12{\omega(v,v')}).
\]
The existence of compatible intertwiners for representations of $H$ is closely related to the fact that the projective Weil representation of the symplectic group $\SP(A)$ can be linearized, cf. \cite{Weil}, \cite{Gerardin}, \cite{MasoudWeil}, \cite{TeruWeil}, \cite{RonnyGeometricWeil}, and especially \cite{RonnyCanonical}.
That is, for gene\-ric choices of $f$,  the automorphism group of $(H,f)$ is isomorphic to $\SP(A)$, and the resulting action of $\SP(A)$ on $\Lambda_f$ lifts to a representation of $\SP(A)$ on $V_f$. This is the Weil representation.

The present work is essentially an extension of this theory to the setting of more general finite nilpotent groups.
\end{rem}

\subsection{Formulas for the intertwiners} \label{formulas}
The representation space $V_{f,\fp}$ is the space of functions $\phi\colon G\ra \bC$ satisfying
$$\phi(pg)=\chi_f(p)\phi(g)$$
for all $g\in G$ and $p\in P$. The representation is given by $(\rho_{f,\fp}(h)\phi)(g)=\phi(gh).$ The averaging operator $\Phi_{\fp_1,\fp_2}$ defined by
\begin{equation}\label{eq:unnorm}
(\Phi_{\fp_1,\fp_2}\phi)(g):= \frac1{\sqrt{|\fp_1||\fp_1\cap\fp_2|}} \sum_{p_1\in P_1} \chi_f^{-1}(p_1) \phi(p_1g)
\end{equation}
is a unitary operator in $\Hom_G(V_{f,\fp_2},V_{f,\fp_1})$.
These intertwiners are not, however, compatible: by Schur's lemma there is a scalar $\alpha(\fp_1,\fp_2,\fp_3)$ of modulus 1 such that
\begin{equation}\label{eq:cocycle}
\Phi_{\fp_1,\fp_2}\circ \Phi_{\fp_2,\fp_3}\circ \Phi_{\fp_3,\fp_1} = \alpha(\fp_1,\fp_2,\fp_3) \Id,
\end{equation}
but in general $\alpha$ does not equal $1$.

Motivated by the theory of the Weil representation, we modify each $\Phi_{\fp_1,\fp_2}$ by a scalar so that the resulting operators are compatible.
As explained in \S\ref{sss:goal}, these scalars depend on a choice of orientations. To a pair $\tfp_1,\tfp_2$ of oriented polarizations, we associate (\S\ref{ss:relativeOrientation}) a \emph{relative orientation} $\theta(\tfp_1,\tfp_2)\in\Fpt/(\Fpt)^2$. Let
\begin{equation}\label{e:defpsi}
\Psi_{\tfp_1,\tfp_2}:=\gamma(1)^{-m(\fp_1,\fp_2)^2} \gamma(\theta(\tfp_1,\tfp_2)) \Phi_{\fp_1,\fp_2},
\end{equation}
where $m(\fp_1,\fp_2)=\log_p |\fp_1|- \log_p|\fp_1\cap \fp_1| - 1$ and
\begin{equation}\label{e:protogamma}
\gamma(a):=\sum_{x\in\Fpt} e^{i \pi ax^2/p}.
\end{equation}

Here is the main result of this paper.
\begin{thm} \label{t:compatible}The operators $\{\Psi_{\tfp_1,\tfp_2}\}$ are compatible; i.e.
\[
\Psi_{\tfp_1,\tfp_2}\circ \Psi_{\tfp_2,\tfp_3}\circ \Psi_{\tfp_3,\tfp_1} =  \Id
\]
for any oriented polarizations $\tfp_1,\tfp_2,\tfp_3$ of $f$.
\end{thm}

\begin{rem}\label{rootindependence} The definitions of $\theta$ and $\gamma$ depend implicitly on the choice of a primitive $p^\th$ root of unity. However, the intertwiner \eqref{e:defpsi} is independent of this choice.
\end{rem}

\subsection{Reduction to the Heisenberg case} \label{ss:ProofReductionProcess}
To prove Theorem \ref{t:compatible}, we first show that the cocycle $\alpha$ of \eqref{eq:cocycle} can be computed in ``Heisenberg terms.''  Our goal now is to formulate this as a precise statement, to be proved in \S\ref{s:canonicalIntertwiners}.
An analysis of $\alpha$ for the Heisenberg group will then yield Theorem \ref{t:compatible}.

Consider $\fg$ as an abelian group equipped with a skew-symmetric pairing $B_f(x,y):=f([x,y]), x,y\in \fg$. Let $\fz\subset\bCt$ be the image of $B_f$.  Equip the abelian group $\bfg:=\fg\oplus \fz$ with the usual Heisenberg commutator: namely, $\fz$ is central and the bracket of $x,y\in \fg\subset \bfg$ equals $B_f(x,y)\in \fz$.

Let $\bf\colon \bfg\ra \bCt$ be the composition $\tfg\onto \fz\inj \bCt$. If $\fp\subseteq \fg$ is a polarization for $f$ then one can show that $\bfp\colon =\fp\oplus \fz$ is a polarization for $\bf$. Using $(\bfg,\bf)$ instead of $(\fg,f)$, we can define $\alpha(\bfp_1,\bfp_2,\bfp_3)$. The main technical result is:

\begin{thm}\label{t:reduction} For all polarizations $\fp_1,\fp_2,\fp_3$ of $f$,
\begin{equation}\label{eq:equalityOfCocycles}
\alpha(\fp_1,\fp_2,\fp_3)=\alpha(\bfp_1,\bfp_2,\bfp_3).
\end{equation}
\end{thm}

To prove this theorem we define, following a suggestion of Drinfeld, the notion of {\it neighboring polarizations} in \S \ref{ss:neighboringPolarizations}. This is the key definition of this paper. We will show that any two polarizations can be connected by a chain of neighbors. Furthermore, the identity \eqref{eq:equalityOfCocycles} is not hard to verify when $\fp_1$ and $\fp_2$ are neighbors (while $\fp_3$ is arbitrary).

We complete the proof of Theorem \ref{t:reduction} in \S\ref{s:canonicalIntertwiners} and the proof of Theorem \ref{t:compatible} in \S\ref{s:compatibleIntertwiners}.

\begin{rem} In the classical setting, the Heisenberg group is associated to the data of a finite dimensional symplectic vector space, as in Remark \ref{r:HeisenbergGroupWeilRep}.
In our setting, the role of the symplectic vector space is played by $(\fg,B_f)$. There are two differences between our setting and the classical one.
First,  $B_f$ may be degenerate, but this is not serious since one can replace $\fg$ by $\fg/\Ker(B_f)$. Second, $\fg$ is not necessarily a vector space. We are, therefore, forced to generalize certain constructions from the setting of vector spaces to finite abelian $p$-groups. We do this in Appendices \ref{s:Det} and \ref{s:WittGroupMaslov}.
\end{rem}

\begin{rem}G. Lion and P. Perrin ~\cite{LP} proved an equality similar to (\ref{eq:equalityOfCocycles}) for unipotent groups over local fields. Using the present notion of neighboring polarizations, one can give a more transparent proof of their important result. \end{rem}

\begin{rem} R. Howe \cite{Howe-Character} also studied automorphisms and representations of $p$-groups by reducing (in a somewhat different way) to the action of a symplectic group on a Heisenberg group.
\end{rem}

\subsection{Organization of the text}
In \S\ref{s:polarizations}, we set out the basic facts about Lie rings, polarizations, orientations, and examine the case of Heisenberg Lie rings. We see that polarizations of a Heisenberg Lie ring are in bijection with Lagrangians in the corresponding symplectic module. We define the notions of neighboring polarizations and relative orientation.
In \S \ref{s:repOfpGroups}, we recall in greater detail the orbit method mentioned in \S\ref{ss:statement}, including the construction due to M. Lazard which defines an equivalence between the relevant categories of nilpotent Lie rings and groups. In \S \ref{s:canonicalIntertwiners} we develop the reduction process of \S\ref{ss:ProofReductionProcess} and prove Theorem \ref{t:reduction}. In \S\ref{s:compatibleIntertwiners} we calculate the cocycle \eqref{eq:cocycle} in the Heisenberg case in terms of the Maslov index (Appendix \ref{s:WittGroupMaslov}), and use the result to prove Theorem \ref{t:compatible}.

In Appendices \ref{s:Det} and \ref{s:WittGroupMaslov}, we explain the necessary generalizations of some constructions in linear algebra to the context of finite abelian $p$-groups. In the literature, these notions are discussed for vector spaces over fields.First, in Appendix \ref{s:Det}, we define and study the determinant functor in two complementary ways. In one approach, following  a suggestion of Drinfeld, we use Deligne's notion of the universal determinant ~\cite[\S4]{DeligneDet} and Quillen's results on $K$-theory ~\cite{Quillen}.  In \S\ref{sec:elementary}, we describe a more elementary approach using filtrations and spectral sequences. In Appendix \ref{s:WittGroupMaslov}, we adapt the theory of the Witt group and the Maslov index to the current setting. Finally, we gather in Appendix \ref{s:proofs} the proofs of some minor results stated in the main text.

\subsection{Conventions}\label{sss:conventions}
Henceforth $p$ is an odd prime;  all the groups, rings, and other structures that appear in this work are implicitly finite and of order a power of $p$. For such a group $A$, we set $\len(A):=\log_p|A|$. If $A$ and $B$ are abelian groups, then $A\otimes B$ denotes the tensor product in the sense of $\bZ$-modules, so a group homomorphism from $A\otimes B$ is the same thing as a bi-additive map from the group $A\times B$.

Finally, the reader should note the following natural but not stan\-dard convention, which causes some (usually simplifying) discrepancies from the literature.   If $e_1,\ldots e_n$ is a basis for a vector space $V$, and $e_1^*,\ldots,e_n^*$ the dual basis for $V^*$, then we identify $\wedge^n(V^*)$ with $(\wedge^nV)^*$ in such a way that
\begin{equation}\label{dualityconvention}(e_1\wedge\cdots\wedge e_n)^*=e_n^*\wedge\cdots\wedge e_1^*.\end{equation}

\subsection{Acknowledgments} We would like to thank V. Drinfeld for introducing us to this project and sharing his notes and insights. In particular, the idea of using neighboring polarizations and results of Deligne and Quillen for defining the determinant are due to him. David Kazhdan and Ronny Hadani helped us understand the importance of this problem and encouraged us to pursue it. Last but not least, we would like to thank Travis Schedler for his interest in this paper and many stimulating conversations.

\section{Lie Rings and Polarizations}\label{s:polarizations}

After laying out the general theory in \S\S\ref{ss:lieRings}--\ref{ss:relativeOrientation}, we will work out the fundamental example (the Heisenberg Lie ring) in  \S \ref{ss:Heisenberg1}.

\subsection{Lie rings}\label{ss:lieRings}
\begin{defe} A {\it Lie ring} is an abelian group $\fg$ equipped with a bracket $[.,.]\colon \fg\otimes \fg\ra \fg$ satisfying the Jacobi identity and the identity $[x,x]=0$.
\end{defe}

\begin{defe} The (Pontryagin) dual of $\fg$ is the abelian group $\fg^*:=\Hom(\fg,\bCt)$. (Note that the bracket of $\fg$ does not enter into the definition of $\fg^*$.)
\end{defe}

\subsection{Polarizations}\label{ss:polarization}
Henceforth $\fg$ denotes a finite {\it nilpotent} Lie ring.

\begin{defe} A {\it polarization} of $f\in \fg^*$ is a Lie subring of $\fp\subseteq \fg$ such that $f([\fp,\fp])=\{1\}$ and $\fp$ is maximal among all {\it subgroups} of $\fg$ with this property.
\end{defe}

From the definition, it is not clear that a functional $f\in \fg^*$ has a polarization. The following theorem of A. Kirillov guarantees the existence of polarizations.

\begin{thm}[\cite{Kirillov}]
\label{t:kirillov}
Every $f\in \fg^*$ has a polarization.
\end{thm}

\subsection{Neighboring polarizations}\label{ss:neighboringPolarizations}
Here is one of the key definitions of this paper.

\begin{defe}Let $\fp_1$ and $\fp_2$ be polarizations for $f\in \fg^*$. We say $\fp_1$ and $\fp_2$ are {\it neighbors} if $[\fp_1,\fp_2]\subseteq \fp_1\cap \fp_2$ (in other words, if $\fp_1$ and $\fp_2$ normalize one another).
\end{defe}

\begin{lemma} \label{l:chainOfPolarizations}
Let $\fp_1$ and $\fp_2$ be polarizations for $f\in \fg^*$. Then there exists a chain of polarizations $\fp_1=\fq_1, \fq_2, \ldots, \fq_m=\fp_2$ such that $\fq_i$ and $\fq_{i+1}$ are neighbors for every $i< m$.
\end{lemma}

The proof is given in \S\ref{pf:chainOfPolarizations}.

\subsection{Relative orientation of oriented polarizations} \label{ss:relativeOrientation}

\begin{defe}\label{d:oriented} An \emph{orientation} of an abelian group $A$ is a non-zero element $o\in\det(A)$  (see Appendix \ref{s:Det}). An {\it oriented abelian group} is a pair $(A,o)$ where $o$ is an orientation of $A$.
\end{defe}

To a pair of oriented polarizations $\tfp_i=(\fp_i,o_i)$, $i=1,2$, we associate an element $\theta(\tfp_1,\tfp_2)\in \Fpt/(\Fpt)^2$, called the \emph{relative orientation}, in the following way. Let $\fl_i:=\fp_i/(\fp_1\cap \fp_2)$. Since $\det$ is an additive functor (\S \ref{ss:additiveFunctors}), we have isomorphisms
 \begin{equation}
\ph_i\colon  \det(\fp_i)\isom \det(\fl_i)\otimes \det(\fp_1\cap \fp_2).
 \end{equation}
 Choose an orientation $\bo$ for $\fp_1\cap \fp_2$, and let $\bo_i$ denote the orientation of $\fl_i$ satisfying $\ph_i(o_i)=\bo_i\otimes \bo$. The pairing
 \[
 B\colon  \fp_1\otimes \fp_2 \ra \bCt, \quad \quad (x,y)\mapsto f([x,y]).
 \]
  induces a perfect\footnote{Suppose $x\in \fp_1$ is such that $B(x,y)=1$ for all $y\in \fp_2$. Then the group generated by $x$ and $\fp_2$ is an isotropic subspace of $\fg$. As $\fp_2$ is a maximal isotropic subgroup, this implies $x\in \fp_2$.} pairing $\fl_1\otimes \fl_2\ra \bCt$, which by Corollary \ref{c:pairingDet}, defines an isomorphism
\begin{equation}\label{lpairing}
b\colon\det(\fl_1)\otimes \det(\fl_2)\isom \Fp.
\end{equation}

\begin{defe} We denote by
$\theta(\tfp_1,\tfp_2)$ the class of $b(\bo_1\otimes\bo_2)$ in $\Fpt/(\Fpt)^2$, and call it the \emph{relative orientation of $\tfp_1$ and $\tfp_2$}.
\end{defe}

\begin{rem} $\theta(\tfp_1,\tfp_2)$  does not depend on the choice of $\bo$ used in its definition. However, as explained in \S\ref{sss:dualityAb}, the isomorphism \eqref{lpairing}, and so $\theta$ itself, depends implicitly on the choice of a primitive $p^\th$ root of unity; see Remark \ref{rootindependence}.
\end{rem}

\begin{rem} Because of our convention \eqref{dualityconvention}, our relative orientation differs from the one common in the literature by a factor of $(-1)^{m(m-1)/2}$, where $m=\len(\fp_1/\fp_1\cap\fp_2)$.
\end{rem}


\subsection{The Heisenberg case} \label{ss:Heisenberg1}

\subsubsection{Heisenberg Lie rings} \label{sss:HeisenbergLieRing}
Let $(A,\omega)$ be a symplectic module; i.e., $A$ is a finite abelian $p$-group and $\omega$ is a symplectic form (see Appendix \ref{ss:modules}). Let $\fz\subset\bCt$ be the image of $\omega$.
The {\it Heisenberg Lie ring} associated to $(A,\omega)$ is the abelian group $\fh:=A\oplus \fz$, equipped with the bracket
\[
[(a,s),(b,t)]=(0,\omega(a,b)),  \quad \quad \quad \forall \, a, b\in A, \quad \forall \, s, t\in \fz.
\]
Note that $[\fh,\fh]=\Z(\fh)=\fz$, where $\Z(\fh)$ denotes the center of $\fh$.

\subsubsection{Polarizations}
\label{sss:PolarizationsForHeisenberg}
Let $f\in \fh^*$ and suppose $\psi\colon \fz\ra \bC^\times$ is the restriction of $f$ to the center $Z(\fh)=\fz$. Note that $\psi \circ \omega$ defines a symplectic pairing on $\bA:=A/\Ker(\psi \circ \omega)$.
Let $\pi\colon A\ra\bA$ be the canonical projection.
Then $$L\mapsto \pi^{-1}(L)\oplus \fz$$
defines a bijection between Lagrangians of $\bA$ and polarizations of $f$.  In particular, if $\psi\circ \omega$ is non-degenerate, then polarizations for $f$ are in bijection with Lagrangians in $A$.

\begin{rem} The pairing $\psi \circ \omega$ is non-degenerate if $\psi$ is injective, and in particular when $f$ is the projection $\fh\to\fz\subset\bCt$.
\end{rem}

\subsubsection{All polarizations for Heisenberg Lie rings are neighbors} \label{sss:neighboringLagrangians}
Suppose $\fp_1$ and $\fp_2$ are polarizations for $f\in \fh^*$. Then
\[[\fp_1,\fp_2]\subseteq [\fh,\fh] =\fz \subseteq \fp_1\cap \fp_2.
\]

\subsubsection{Oriented polarizations}\label{sss:relativeOrientLagrangians}
Let $\tL_1,\tL_2$ be oriented Lagrangians in $A$.
We can define the relative orientation $\theta(\tL_1,\tL_2)$ by modifying the procedure of \S\ref{ss:relativeOrientation}
in an obvious way.

Suppose that $\psi\circ\omega$ is non-degenerate.
Let $\fp_i:=L_i \oplus \fz$. According to \S \ref{sss:PolarizationsForHeisenberg}, each $\fp_i$ is a polarization for $f$.
Choosing an  orientation $o\in \det(\fz)$, we obtain an orientation $o_i\in \det(\fp_i)=\det(L_i)\otimes \det(\fz)$. It is easy to show that
$\theta(\tfp_1,\tfp_2)=\theta(\tL_1,\tL_2)$.

\begin{rem} The determinant $\det(\fz)$ is equal to the group of $p^\th$ roots of unity in $\bC$. So an orientation of $\fz$ is just a primitive $p^\th$ root of unity.
\end{rem}

\section{The Orbit Method}\label{s:repOfpGroups}

Here we describe in more detail the orbit method of \S\ref{ss:statement}, following the approach of ~\cite[\S2]{DB}. The example of the Heisenberg group is worked out in \S\ref{ss:Heisenberg2}.

\subsection{The Lie ring of a finite $p$-group}\label{ss:LieRingOfGroups}
Let $\nilp_p$ denote the category of nilpotent Lie rings of order a power of $p$ and nilpotence class less than $p$. To every $\fg\in \nilp_p$ we associate a finite $p$-group $G=\Exp(\fg)$ as follows. $\Exp(\fg)$ has the same underlying set as $\fg$ and the group operation is defined by
\begin{equation}
x*y:=\sum_{i<p} \CH_i(x,y)      \quad \quad  \quad       x,y\in \fg,
\end{equation}
where $\CH_i$ is the homogenous component of degree $i$ of the Campbell-Hausdorff
series
$$\begin{aligned} \CH(x,y)&:=\log(\exp(x)\exp(y))\\ &=x+y+\frac12[x,y]+\mbox{(higher degree terms)}.\end{aligned}$$ It is well-known that $\CH_i$ is a polynomial with coefficients in $\bZ[\frac{1}{i!}]$. Hence the group operation is well-defined. It is easy to check that if $\varphi\colon \fg\ra \fg'$ is a morphism of Lie rings, then $\varphi$, viewed as a map $\Exp(\fg)\ra \Exp(\fg')$, is a group homomorphism. Let $\Nilp_p$ denote the category of finite $p$-groups of nilpotence class less than $p$.

\begin{thm}[\cite{Lazard}]
$\Exp\colon \nilp_p\ra \Nilp_p$ is an equivalence of categories.
\end{thm}
We will denote the inverse of $\Exp$ by $\Log$. If $G\in \Nilp_p$, we refer to $\Log(G)$ as the Lie ring associated to $G$. We denote the identity maps between $\fg$ and $G$ by $\log\colon G\ra \fg$ and $\exp\colon  \fg \ra G$.

\subsection{Coadjoint orbits} \label{ss:coadjoingOrbits}
Let $G\in \Nilp_p$ and let $\fg:=\Log(G)$. $G$ acts on $\fg$ by conjugation. This is known as the {\it adjoint action}. The corresponding right action of $G$ on $\fg^*:=\Hom(\fg,\bCt)$ is known as the {\it coadjoint action}. The orbits of the action of $G$ on $\fg^*$ are known as the {\it coadjoint orbits}.

\subsection{The orbit method}\label{s:orbitMethod}
Fix $G\in \Nilp_p$, $\fg:=\Log(G)$ and $f\in \fg^*$ and let $\Omega \subseteq \fg^*$ denote the coadjoint orbit containing $f$. Let $\fp$ be a polarization for $f$ and let $P:=\Exp(\fp)\subseteq G$. Then $\chi_f\colon P\ra \bCt$ defined by
$$\chi_f(x)=f(\log x)$$ is a one-dimensional representation of $P$. Let $\rho_{f,\fp}:=\ind_P^G \chi_f$.
As in \S\ref{formulas}, we write $V_{f,\fp}$ for the representation space of $\rho_{f,\fp}$.

\begin{thm} [\cite{Kirillov}, \cite{DB}, \S 2]\label{t:orbitMethod}
The representation $\rho_{f,\fp}$ is irreducible and its character is given by \[
g\mapsto \frac{1}{\sqrt{|\Omega|}} \sum_{f\in \Omega} f(\log g), \quad \quad \forall \, g\in G.
\]
For each coadjoint orbit $\Omega$, choose $f\in \Omega$ and a polarization $\fp$ of $f$. Then the map $\Omega \mapsto \rho_\Omega:=\rho_{f,\fp}$ is a bijection between the set of coadjoint orbits and the set of isomorphism classes of irreducible representations of $G$.
\end{thm}

\subsection{Example: Heisenberg groups}\label{ss:Heisenberg2}

Suppose $\fh$ is the Heisenberg Lie ring associated to a symplectic module $(A,\omega)$  as in \S\ref{ss:Heisenberg1}. Then the underlying set of the Heisenberg group  $H:=\Exp(\fh)$ is $A \times \fz$ and the multiplication is given by
\[
(a,t)*(b,s)=(a+b, t\cdot s\cdot\omega(a,b)^{1/2})
\]
for all $a,b\in A$ and $s,t\in \fz$. (Note that $\fz$ is written multiplicatively, and $\omega(a,b)$ has a unique square-root in $\fz$.)

\subsubsection{Coadjoint orbits} \label{sss:heisenbergcoadjoint}
The adjoint action $H$ on $\fh$ is easy to compute: given $(a,t)\in H$, we have
\[
\Ad(a,t)\colon  \fh\ra \fh, \quad \quad (b,s)\mapsto (b, s\cdot\omega(a,b)).
\]
The coadjoint action of $H$ on $f\in \fh^*$ is given by
\[
[\Ad^*(a,t)f](b,s):=f([\Ad(a,t)](b,s))=f(b,s\cdot\omega(a,b)).
\]
Let $\psi$ denote the restriction of $f$ to $Z(\fh)=\fz$. Then
the orbit containing $f$ has size $|A/\Ker(\psi \circ \omega)|$.

\subsubsection{The Stone-von Neumann Theorem} \label{sss:SvN} Theorem \ref{t:orbitMethod} has the following special case, known traditionally as the Stone-von Neumann Theorem. Choose a non-trivial character $\psi\colon \fz\to \bCt$, and define $f\in \fh^*$ by $f(v,t)=\psi(t)$.
 According to \S\ref{sss:PolarizationsForHeisenberg},  polarizations of $f$ are of the form $\fp:=L\oplus \fz$ where $L$ is a Lagrangian in $A$. The orbit method states that the corresponding representation $\rho_{f,\fp}:=\ind_P^H \chi_f$ is irreducible, and its isomorphism class does not depend on $L$. Moreover, this is the unique isomorphism class of representations with central character $\psi$.

\section{The Reduction Process}
\label{s:canonicalIntertwiners}

\subsection{Reduction of Lie rings to the Heisenberg case}
\label{ss:reduction}
Let $\fg$ be a  Lie ring. For every $f\in \fg^*$, let $B_f\colon  \fg\otimes \fg\ra \bCt$ denote the skew pairing
\begin{equation}
B_f(x,y)=f([x,y]), \quad \quad \forall \, x,y\in \fg.
\end{equation}
Let $\omega$ be the induced symplectic form on
$\fg/\Ker(B_f).$

\begin{defe}
We will denote by  $\bfg$  the Heisenberg Lie ring associated to the pair $(\fg/\Ker(B_f),\omega)$ as in \S\ref{ss:Heisenberg1}, so  $\bfg$ equals $\fg/\Ker(B_f)\oplus \fz$ as an abelian group.  If $\fp\subseteq \fg$ is a subgroup containing $\Ker(B_f)$, we denote by $\bfp\subseteq \bfg$ the subgroup $\fp/\Ker(B_f)\oplus \fz$. Finally, let $\bf\in (\bfg)^*$ denote the composition
\[\bf\colon \bfg\onto \fz \subset \bCt. \]

\end{defe}

\subsubsection{Polarizations in the reduction process} \label{sss:PolarizationinReduction}

If $\fp$ is a polarization for $f$ then $\bfp$ is a polarization for $\bf$.  By the discussion of \S\ref{sss:PolarizationsForHeisenberg}, $\bfp=L\oplus \fz$ for some Lagrangian $L\subset\bfg$.

\begin{rem}\label{r:sizeOfPolarization} The fact that $\fp/\Ker\,B_f$ is a Lagrangian in $\fg/\Ker\,B_f$ shows that $|\fp|=\sqrt{|\fg||\Ker\,B_f|}$.
\end{rem}

\subsubsection{Invariance of relative orientation under the reduction process} \label{sss:relativeOrientReduction}
Let $\tfp_i=(\fp_i,o_i)$, $i=1,2$, be oriented polarizations for $f$.
Choose (once and for all) orientations of  $\Ker(B_f)$ and $\fz$.
Then the isomorphisms
$$\begin{aligned}
 \det(\fp)&\cong \det(\fp/\Ker(B_f))\otimes \det(\Ker(B_f))\\
 \det(\fp^f) & \cong \det(\fp/\Ker(B_f)) \otimes \det(\fz).
 \end{aligned}
 $$
determine orientations $o_i^f$ of $\fp_i^f.$
Let $\tfp_i^f:=(\fp_i^f,o_i^f)$. It is easy to show that $\theta (\tfp_1^f,\tfp_2^f)=\theta(\tfp_1,\tfp_2),$ independent of the  choices made.


\subsection{The cocycle of three polarizations}\label{ss:cocycle}
Let $\fg$ be a Lie ring of nilpotence class less than $p$, and $G:=\Exp(\fg)$ the corresponding finite $p$-group, as in \S\ref{ss:LieRingOfGroups}. Fix $f\in \fg^*$.
For each polarization $\fp$ of $f$,  let $(V_{f,\fp},\rho_{f,\fp})$ denote the corresponding irreducible representation of $G$, as in \S\ref{s:orbitMethod}. For each pair $\fp_1,\fp_2$ of polarizations we have the intertwiner $\Phi_{\fp_1,\fp_2}\in\Hom_G(V_{f,\fp_2},V_{f,\fp_1})$
defined by \eqref{eq:unnorm}. For each triple $\fp_1,\fp_2,\fp_3$ we have  the unitary scalar
 $\alpha(\fp_1,\fp_2,\fp_3)$ defined by \eqref{eq:cocycle}.
Our goal is to prove Theorem \ref{t:reduction}, i.e. that
$$\alpha(\fp_1,\fp_2,\fp_3)=\alpha(\fp_1^f,\fp_2^f,\fp_3^f).$$
We proceed in three steps.

\underline{Step 1: Basic formula}\label{sss:basic}  First let us give a clean formula for $\alpha(\fp_1,\fp_2,\fp_3)$. Let
$$S=S(\fp_3,\fp_2,\fp_1):=\{(p_3,p_2,p_1)\mid p_i\in P_i, p_3p_2p_1=1\}.$$
Simple manipulations with \eqref{eq:unnorm}, explained in \S\ref{pf:symfla}, yield the following lemma.
\begin{lemma}\label{p:symfla} For any polarizations $\fp_1,\fp_2,\fp_3$ of $f$,
\begin{equation}\label{eq:form}\alpha(\fp_1,\fp_2,\fp_3)=N(\fp_1,\fp_2,\fp_3)\sum_{(p_3,p_2,p_1) \in S}(\chi_f(p_3) \chi_f(p_2)\chi_f(p_1))^{-1}\end{equation}
where $$N(\fp_1,\fp_2,\fp_3):=
\sqrt{\frac{|\fg||\Ker\, B_f|}{|\fp_1||\fp_2||\fp_3||\fp_1\cap\fp_2||\fp_2\cap\fp_3||\fp_3\cap\fp_1|}}.$$
\end{lemma}

\underline{Step 2: Reduction to neighboring polarizations}  By Lemma \ref{l:chainOfPolarizations}, we can
find a chain $\fp_1=\fq_1, \ldots,\fq_m=\fp_2$ of polarizations such that consecutive ones are neighbors. Using the fact that
$\Phi_{\fq_i,\fq_{i+1}}=\Phi_{\fq_{i+1},\fq_i}^{-1}$ we can write
$$\alpha(\fp_1,\fp_2,\fp_3)=\prod_{i=2}^{m-1}\alpha(\fq_i,\fq_{i+1},\fp_1)\cdot  \prod_{i=1}^{m-1}\alpha(\fq_{i+1},\fq_{i},\fp_3).$$
In each factor, the first two polarizations are neighbors. So we may
as well assume that  $\fp_1$ and $\fp_2$ are neighbors themselves.

\underline{Step 3: Computation for neighboring polarizations}
Suppose that $\fp_1$ and $\fp_2$ are neighbors.  Let
$$\fS=\{(p_3,p_2,p_1)\mid p_i\in \fp_i, p_3+p_2+p_1=0\},$$ and let $\fR$ be the relation on $\fS$ given by
\begin{equation*}\label{rel2}(p_3,p_2,p_1)\sim(a_{13}+p_3-a_{32},a_{32}+p_2-a_{21},a_{21}+p_1-a_{13})\end{equation*}
for all selections of $a_{ij}\in \fp_i\cap \fp_j$.

\begin{lemma}\label{l:neighborfla}
When $\fp_1,\fp_2$ are neighbors, and $\fp_3$ arbitrary,
\begin{equation}\label{eq:cocycleNeighbors}\alpha(\fp_1,\fp_2,\fp_3)=N'(\fp_1,\fp_2,\fp_3) \sum_{(p_3,p_2,p_1)\in \fS/\fR} f(\frac{1}{2}[p_2,p_1])\end{equation}
where
$$N'(\fp_1,\fp_2,\fp_3):=\sqrt{\frac{|\fg||\fp_1\cap\fp_2||\fp_2\cap\fp_3||\fp_3\cap\fp_1|}{|\fp_1\cap\fp_2\cap\fp_3|^2|\fp_1||\fp_2||\fp_3|}}\sqrt{|\Ker\, B_f|}.$$
\end{lemma}

The proof is given in \S\ref{pf:neighborfla}.
 Note that the right-hand side of \eqref{eq:cocycleNeighbors} remains the same if $\fg,f,\fp_i$ are replaced by $\bfg,\bf,\bfp_i$.  We thus conclude the proof of Theorem \ref{t:reduction}.

\section{Compatible Intertwiners}\label{s:compatibleIntertwiners}

In this section we restate and prove Theorem \ref{t:compatible}.  First we calculate the cocycle $\alpha$ of \eqref{eq:cocycle}  in the case of a Heisenberg group. Then we will apply Theorem \ref{t:reduction}.

\subsection{The Heisenberg case}\label{ss:Heisenberg3} Let $(A,\omega)$ be a symplectic module and $\fh$ the corresponding Heisenberg Lie ring, as in \S\ref{ss:Heisenberg1}. Fix $f\in \fh^*$ and let $\psi$ denote the restriction of $f$ to the center $\fz$ of $\fh$. Assume $\psi$ is injective. By
\S\ref{sss:PolarizationsForHeisenberg}, the map $L\mapsto L\oplus \fz$ is a bijection between Lagrangians in $A$ and polarizations of $f$. Let $\tau(L_1,L_2,L_3)$ denote the Maslov index of the triple $(L_1,L_2,L_3)$ (see \S\ref{ss:Maslov}) and let $\gamma$ denote Weil's gamma index (\S \ref{ss:GammaIndex}).

\begin{prop} \label{p:cocycleIsMaslov} For any Lagrangians $L_1,L_2,L_3$ in $A$,
$$\alpha(L_1\oplus\fz,L_2\oplus\fz,L_3\oplus\fz)=\gamma(\tau(L_1,L_2,L_3)).$$
\end{prop}

The proof is given in \S\ref{pf:cocycleIsMaslov}.

\subsection{The Compatibility Theorem}

Fix $\fg$, $G$, and $f$ as in \S  \ref{ss:cocycle}.
For any oriented polarizations $\tfp_1,\tfp_2$ of $f$,  define
\begin{equation} \label{eq:canonicalIntertwiners}
\beta(\tfp_1,\tfp_2):=\gamma(1)^{-m(\fp_1,\fp_2)^2} \gamma( \theta(\tfp_1,\tfp_2) ),
\end{equation}
where $m(\fp_1,\fp_2):=\len(\fp_1)-\len(\fp_1\cap \fp_2)-1$,  and $\theta(\tfp_1,\tfp_2)$ is the relative orientation of $\tfp_1$ and $\tfp_2$ defined in \S \ref{ss:relativeOrientation}.\footnote{Note that $\gamma(a)$ depends only on the class of $a$ in $\Fpt/(\Fpt)^2$. In particular, the formula for $\beta$ is well-defined.}  Our goal is to prove Theorem \ref{t:compatible}, which amounts to the identity
\begin{equation}
\label{eq:cocycleIscoboundary}
\alpha(\fp_1,\fp_2,\fp_3)=\beta(\tfp_1,\tfp_2)\beta(\tfp_2,\tfp_3)\beta(\tfp_3,\tfp_1)
\end{equation}
for any triple $\tfp_1,\tfp_2,\tfp_3$ of oriented polarizations of $f$.

\subsubsection{}
Theorem \ref{t:reduction} states that the left-hand side of (\ref{eq:cocycleIscoboundary}) is unchanged if we replace $\fp_i$ by $\fp_i^f$. On the other hand, by the discussion of \S \ref{sss:relativeOrientReduction}, the right-hand size is unchanged if we replace $\tfp_i$ by $\tfp_i^f$. In view of Proposition \ref{p:cocycleIsMaslov},
formula \eqref{eq:cocycleIscoboundary} is a consequence of the following proposition about the Maslov index.

\begin{prop} \label{p:intertwinersforHeisenberg}
Suppose $\tL_1,\tL_2,\tL_3$ are oriented Lagrangians in a symplectic module. Then,
\[
\beta(\tL_1,\tL_2)\beta(\tL_2,\tL_3)\beta(\tL_3,\tL_1)=\gamma(\tau(L_1,L_2,L_3)).
\]
\end{prop}
The proof is given in \S\ref{pf:intertwinersforHeisenberg}.

\appendix

\section{Determinants of Finite Abelian Groups}\label{s:Det}

The determinant of a finite abelian $p$-group $A$ will be defined as a one-dimensional vector space over $\Fp$ (and the determinant of a general finite abelian group could be defined as the sequence of determinants of its Sylow subgroups). In the case when $A$ is a vector space over $\Fp$, its determinant is simply the top exterior power of $A$.

In \S \S \ref{ss:KTheoryReview}--\ref{ss:det} we develop this theory in the natural setting of the $K$-theory of exact categories, following Quillen \cite{Quillen} and Deligne \cite{DeligneDet}. In \S\ref{sec:elementary} we give a second, more elementary treatment in terms of filtrations and spectral sequences. For an approach based on exterior powers see \cite[\S2]{Reichstein2002}.

\subsection{Recollections on $K$-theory} \label{ss:KTheoryReview}

\subsubsection{The $K$-theory space}\label{sss:KTheorySpace}
Let $\cM$ be an exact category. To define the algebraic $K$-theory of $\cM$, Quillen ~\cite{Quillen} defined a category $\Q \cM$ which has the same objects as $\cM$. (We do not need the precise description of morphisms in this category.) Let $\BQM$ denote the geometric realization of $\Q \cM$. Direct sums in $\cM$ induce operations $\QM\times \QM\ra \QM$ and $\BQM\times \BQM\ra \BQM$ endowing $\BQM$ with a structure of a commutative $H$-space.

Fix a zero object $0\in \cM$. The $K$-theory space $\K \cM$ of $\cM$ is the based loop space of $\BQM$. One sets $\K_i(\cM):=\pi_i(\K \cM,0)$.

\subsubsection{The category of exact sequences} \label{sss:categoryOfExactSeq}
Let $\cE$ denote the category of exact sequences of the exact category $\cM$. For an exact sequence
$$\Sigma\colon  M'\inj M\onto M''$$
 let $s\Sigma$, $t\Sigma$, $q\Sigma$ denote the sub-, total, and quotient objects of $M$. A sequence in $\cE$ is called \emph{exact} if it gives rise to three exact sequences in $\cM$ on applying $s$, $t$, $q$. The category $\cE$ is thus endowed with the structure of an exact category.

\begin{thm}[\cite{Quillen}, Theorem 2]\label{t:categoryOfSES}
The functor $(s,q)\colon  \Q(\cE)\ra \Q(\cM)\times \Q(\cM)$ is a homotopy equivalence.
\end{thm}

\subsubsection{Devissage}\label{sss:devissage}
Suppose $\cM$ is an Artinian abelian category. Let $\cS$ denote the full subcategory of semisimple objects.
\begin{thm}[\cite{Quillen}, Theorem 4]\ \label{t:KTheoryOfArtinians}
The natural functor $\Q(\cS)\ra \Q(\cM)$ is a homotopy equivalence.
\end{thm}

\subsection{Picard groupoids}\label{ss:picardGroupoids}
A {\it Picard groupoid} is a  monoidal category in which all objects and morphisms are invertible. A 1-morphism of Picard groupoids (or a {\it{Picard functor}}) is a monoidal functor between the corresponding monoidal categories. A Picard groupoid is commutative if the corresponding monoidal category is symmetric.

Let $\cP$ be a Picard groupoid. We denote the set of isomorphism classes of $\cP$ by $\pi_0(\cP)$. Let $\pi_1(\cM)$ denote the group of automorphisms of the unit object of $\cM$. One can show that $\pi_1(\cM)$ is abelian.

\subsubsection{Super lines} \label{sss:superLines}
Let $R$ be a commutative ring. The commutative Picard groupoid $\Pic_R$ of super lines on $R$ is defined as follows.
The objects of $\Pic_R$ are pairs $(L,n)$ where $L$ is a rank-one free module over $R$ and $n$ is an integer. A morphism $(L,m)\to (L',m')$ exists only if $m=m'$, in which case it is any isomorphism $L\to L'$ of modules. The monoidal product is
\[
(L,n)\otimes (L',n'):= (L\otimes L', n+n').
\]
The inverse of $(L,n)$ is ($L^\vee, -n)$, where $L^\vee:=\Hom(L,R)$ is the dual line.
Finally, the commutativity constraint $(L,n)\otimes (L',n') \ra (L',n')\otimes (L,n)$ is defined using ``the Koszul sign rule'':
 \begin{equation}\label{eq:koszul}
 a\otimes a'\mapsto (-1)^{nn'} a'\otimes a, \quad \quad \mbox{for $a\in L, b\in L'$}.
 \end{equation}

 \begin{rem} \label{r:dualityPic}
 (Duality) The assignment $(L,n)\mapsto (L^\vee, n)$ defines a Picard functor
$\DPic: \Pic_F\ra \Pic_F^o$, which is an equivalence of Picard categories.
 \end{rem}

 \begin{rem} \label{r:piOfSuperLines}
 Note that $\pi_0(\Pic_R)=\bZ$ and $\pi_1(\Pic_R)=R^\times$.
 \end{rem}

\subsubsection{The fundamental groupoid of the $K$-theory space}\label{sss:fundamentalGroupoid}
 Let $\cM$ be an exact category. Let $\VM$ denote the fundamental groupoid of the $K$-theory space of $\cM$: the objects of $\VM$ are loops $\gamma$ in $\BQM$ based at $0$. A morphism $\gamma \ra \gamma'$ is a homotopy class of based homotopies $\gamma \ra \gamma'$.  Following Deligne \cite{DeligneDet}, we call $\VM$ the category of virtual objects of $\cM$.
 The composition of loops defines a multiplication on $\PKM$ endowing it with a structure of a Picard groupoid.
The usual proof that the fundamental group of a group is commutative shows that $\PKM$ is naturally commutative as a Picard groupoid, being the fundamental groupoid of an $H$-space (\S\ref{sss:KTheorySpace}).

\subsection{Additive functors}\label{ss:additiveFunctors}
Let $\cM$ be an exact category. Let $(\cM,\is)$ be the groupoid whose objects are the same as $\cM$ and whose arrows are isomorphisms in $\cM$. Let $\cP$ be a Picard groupoid.

\begin{defe}\label{def:additivefunctor} A {\it additive functor} $\cM\ra \cP$ is the data of
\begin{enumerate}
\item[(a)] a functor $[\quad]\colon  \cMo \ra \cP$,
\item[(b)] for every exact sequence $\Sigma$, an isomorphism
\[
\{\Sigma\}\colon  [M]\isom [M']+[M''],
\]
functorial for isomorphisms of exact sequences, and
\item[(c)] for every zero object of $\cM$, an isomorphism $[0]\isom 0$,
\end{enumerate}
satisfying the following axioms:
\begin{enumerate}
\item If $\ph\colon M\ra N$ is an isomorphism and $\Sigma$ is the exact sequence $0\ra M \ra N$ (resp. $M\ra N\ra 0$) Then $[\ph]$ (resp. $[\ph]^{-1}$) equals the composition
\[
\quad \quad \quad \quad [M]\rar{\{\Sigma\}} [0] + [N] \rar{} [N]
\]
\[
\mathrm{(resp. } \quad \quad [N] \rar{ \{\Sigma \} } [M] + [0] \rar{  } [M] \mathrm{).}
\]

\item If $0\subseteq M\subseteq N\subseteq P$ is an admissible filtration in $\cM$, the diagram of isomorphisms (coming from (b))
\[
\xymatrix{
[P] \ar[d] \ar[r]& [M]+[P/M]\ar[d] \\
[N]+[P/N] \ar[r] & [M]+[N/M]+[P/N]
}
\]
commutes.

\end{enumerate}
\end{defe}

\begin{rem} An additive functor to a commutative Picard groupoid is \emph{compatible with the commutativity} if for every $M',M''\in\cM$, the diagram
\[
\xymatrix{
 [M']+[M''] \ar[dr]_{\{\Sigma\}} \ar[rr] & & [M'']+[M']\ar[dl]^{\{\Sigma' \}}\\
 &                      [M'\oplus M''] &
 }
\]
commutes. Here the horizontal arrow is the commutativity isomorphism in $\cP$, $\Sigma$ is the exact sequence $M'\inj M'\oplus M''\onto M''$, and $\Sigma'$ is the exact sequence $M'' \inj M'\oplus M''\onto M'$. \end{rem}

\subsubsection{Universal additive functor}\label{sss:universalFunctor}
Let $\cM$ be an exact category and let $\cE$ denote its category of exact sequences.  We define an additive functor $\Delta\colon \cM\to \PKM$ compatible with the commutativity. Given $M\in \cM$, the morphisms $0\inj M$ and $M\onto 0$ define a loop $\Delta(M)\in \BQ\cM$ based at $0$, and hence an object of $\PKM$. Let $M' \inj M \onto M''$ be an exact sequence in $\cM$. The homotopy equivalence $B\cE \to \BQM\times \BQM$ (Theorem \ref{t:categoryOfSES}) defines an isomorphism $\Delta(M) \to  \Delta(M')+\Delta(M'')$ in $\PKM$ making $\Delta$ into an additive functor.

\begin{thm}[\cite{DeligneDet}, \S 4] \label{t:UniversalDet}
Let $\cP$ be a Picard groupoid. Then composition with $\Delta$ defines an equivalence of groupoids between  Picard functors $\PKM \ra \cP$ and additive functors $\cM\ra \cP$.
\end{thm}

\subsection{Determinants}\label{ss:det}

\subsubsection{Determinants of vector spaces} \label{sss:determinant}
Let $F$ be a field, and let $\Vect_F$ denote the abelian category of finite dimensional vector spaces over $F$. The determinant $\Det\colon \Vect_F\ra \Pic_F$ assigns to every vector space $V$ the super line $(\det(V),n)$, where $n=\dim(V)$, and $\det(V)=\wedge^nV$ denotes the $n^\th$ exterior power of $V$. One checks that the determinant is an additive functor which is compatible with commutativity.  This is the reason for choosing the Koszul rule of signs (\ref{eq:koszul}).

\begin{rem} \label{r:Pic=PiKVect}
According to Theorem \ref{t:UniversalDet}, $\Det\colon  \Vect_F\ra \Pic_F$ factors through a Picard functor $\V(\Vect_F)\ra \Pic_F$. The latter functor induces isomorphisms on $\pi_0$ and $\pi_1$, hence it is an equivalence of Picard groupoids.
\end{rem}

\subsubsection{Compatibility with duality in $\Vect_F$} \label{sss:dualityVect} Recall that the duality functor
\[
\fD_{\Vect_F}\colon  \Vect_F \ra \Vect_F^o,\quad \quad V\mapsto {V^\vee}:=\Hom(V,F),
 \]
 is an equivalence of exact categories.
Our convention \eqref{dualityconvention} defines a natural isomorphism
\begin{equation}\label{dualconvention} \Det\circ\fD_{\Vect_F} \cong \DPic\circ \Det
\end{equation}
where $\DPic$ was defined in Remark \ref{r:dualityPic}.

 \subsubsection{Determinants of abelian p-groups}\label{sss:detAb}
Let $\Abp$ denote the category of finite abelian $p$-group. Note that $\Vect_{\Fp}$ is the full subcategory of semisimple objects in $\Abp$. By Theorem \ref{t:KTheoryOfArtinians}, the natural functor $\Q \Abp \ra \Q \Vect_{\Fp}$ is a homotopy equivalence. In particular, we have an equivalence of Picard groupoids $\V (\Abp) \eqv \V (\Vect_{\Fp})$.

\begin{defe} The determinant functor $\Det\colon  \Abp\ra \Pic_{\Fp}$ is the composite
 \[
 \Abp \rar{\Delta} \V (\Abp) \eqv \V(\Vect_{\Fp}) \eqv \Pic_{\Fp}.
 \]
We define $\det(A)$ to be the one-dimensional vector space such that $\Det(A)=(\det(A),\len(A))$.
\end{defe}

\subsubsection{Compatibility with duality in $\Abp$} \label{sss:dualityAb}
The Pontryagin duality functor
 \[\fD_p\colon  \Abp\ra \Abp^o \quad \quad A\mapsto A^*:=\Hom(A,\bCt)
 \]
 is an equivalence of categories. Let $i$ be the inclusion $\Vectp\inj \Abp$.  The choice of a non-trivial character  $\psi_0\colon  \Fp\ra \bCt$ determines
a natural isomorphism $\alpha\colon\fD_{\Vectp} \ra \fD_p\circ i$. (One should use the same $\psi_0$ as in Example \ref{e:angles}.)
This in turn determines a natural isomorphism
 $\Det\circ\fD_p = \DPic\circ \Det$ as follows.
\begin{equation*}
\xymatrix{
 \Abp \ar[r]^{\Delta}\ar[d]^{\fD_p} &  \V (\Abp) \ar[r] \ar[d]^{\V(\fD_p)}&  \V(\Vect_{\Fp}) \ar[r] \ar[d]^{\V(\fD_{\Vect_p})} & \Pic_{\Fp} \ar[d]^{\DPic}.
 \\
 \Abp^o \ar[r]^{\Delta} &  \V (\Abp^o) \ar[r] &  \V(\Vect_{\Fp}^o) \ar[r] & \Pic_{\Fp}^o.
}
\end{equation*}
Here the composition along the top and down is
 $\DPic\circ \Det$, while the composition down and along the bottom is
 $\Det\circ\fD_p$. The left-hand square commutes canonically by the functoriality of $\V$; the middle square commutes up to the natural isomorphism $\alpha$; and the right-hand square commutes up to the natural isomorphism \eqref{dualconvention}.

\begin{cor} \label{c:pairingDet}
Let $\beta\colon  A\otimes B\ra \bCt$ be a non-degenerate pairing. Then $\beta$ defines an isomorphism $\det(\beta) \colon  \det(A)\otimes \det(B) \isom \Fp$.
\end{cor}

\begin{rem} As we have defined it, $\det(\beta)$ depends on the choice of the character $\psi_0\colon\Fp\to\bCt$. In fact, writing $\mu_p$ for the group of $p^\th$ roots of unity, one gets an isomorphism
$\det(A)\otimes \det(B) \isom \mu_p^{\otimes \len(A)}$
independent of $\psi_0$.
\end{rem}

\subsection{Elementary description of the determinant}\label{sec:elementary}

In this section we give an elementary definition of the determinant functor $\Abp\to\Pic_{\Fp}$.

\begin{defe} Given $A\in \Abp$ define the \emph{canonical filtration} $A_\bullet$ of $A$ by $A_n:=p^nA$. Each $\gr_n(A):=A_n/A_{n+1}$ is a vector space over $\Fp$. Define
$$\det(A):=\det(\gr_0(A))\otimes\det(\gr_1(A))\otimes \cdots\otimes\det(\gr_N(A))\in \Pic_{\Fp}$$
for any $N$ such that $p^NA=0$.
\end{defe}

\subsubsection{} Let us explain how to upgrade this function $\det$ to an additive functor $\Abp\to\Pic_{\Fp}$.
Any isomorphism $A\to B$ in $\Abp$ respects the canonical filtration, so $\det$ defines a functor $[\quad]:=(\det,\len)$ as in Definition \ref{def:additivefunctor}(a). The isomorphism $\det(0)\to \Fp$ of Definition \ref{def:additivefunctor}(c) is also automatic. It remains, as in Definition \ref{def:additivefunctor}(b), to construct an isomorphism $\{\Sigma\}\colon \det(K^0)-\det(K^1)+\det(K^2)\to 0$ from any exact sequence
\begin{equation}
\xymatrix{ 0 \ar[r] & K^0 \ar[r] & K^1 \ar[r] & {K^2} \ar[r] & 0.}
\end{equation}
Suppose first that $V^\bullet:=(V^0\to\cdots\to V^n)$ is a finite complex of vector spaces. Then the additive functor $\Det\colon \Vect_{\Fp}\to\Pic_{\Fp}$ induces a canonical isomorphism $\{V^\bullet\}$ from
$$\det(V^\bullet):=\det(V^0)-\det(V^1)+\cdots+(-1)^n\det(V^n)$$
to
$$\det(H^0(V^\bullet))-\det(H^1(V^\bullet))+\cdots+(-1)^n\det(H^n(V^\bullet)).$$
(Condition (2) of \ref{def:additivefunctor} ensures these isomorphisms are well-defined.)
Now a standard spectral sequence construction\footnote{See ~\cite{LA}, Chapter XX, Proposition 9.3. Notationally, his $F^iK^\bullet$ is our $K^\bullet_i$, and we have put $E^n_m:=\bigoplus_{p+q=n}E_m^{p,q}.$}  gives a sequence $E^\bullet_0,$ $E^\bullet_1,\ldots,E^\bullet_N$ of complexes of vector spaces such that
$$E_0^n=\bigoplus_{i=0}^{\infty}\gr_i(K^n), \qquad E_m^n=H^n(E_{m-1}^\bullet)\mbox{ for $m>0$,} \qquad E_N^\bullet=0.$$ The various isomorphisms $\{E_m^\bullet\}$ therefore compose to give an isomorphism
$$\alpha_1\colon \det(E_0^\bullet)\to \det(E_2^\bullet) \to \cdots\to \det(E_N^\bullet)=0.$$
On the other hand, there is also a canonical isomorphism
$$\alpha_2\colon \det(K^0)-\det(K^1)+\det(K^2)\to \det(E_0^\bullet)$$
induced by the canonical isomorphisms $\det(K^n)\cong \det(E_0^n)$.  To complete the construction of the additive functor $\det$, set $\{\Sigma\}:=\alpha_1\circ\alpha_2$.

\subsubsection{} The conditions $(1,2)$ of Definition \ref{def:additivefunctor} can be reduced to the corresponding properties of the determinant functor on vector spaces. The compatibility of $\det$ with duality, as in Corollary \ref{c:pairingDet}, is more subtle, and we leave it to the reader.

\section{Witt Group, Maslov Index, Gamma Index}\label{s:WittGroupMaslov}

In this section we study finite abelian $p$-groups equipped with non-degenerate symmetric forms; we call these \emph{quadratic modules}.
In \S\ref{ss:modules} we set out some general conventions, and in \S\ref{sss:discriminant} define the discriminant of a quadratic module as an element of $\Fptt$. In \S\ref{ss:wittGroup} we define the Witt group $\bW$ of quadratic modules. The analogue for quadratic forms on vector spaces is very well known, cf. \cite{Lam}. In \S\ref{ss:Maslov} we study the Maslov index, an invariant associated to a collection of Lagrangian subspaces of a symplectic module. Again, the analogue for vector spaces is well known; we follow the development in \cite{TeruMaslov}. Finally, in \S\ref{ss:GammaIndex}, we define a character $\gamma$ of $\bW$, following \cite{Weil}.

\subsection{Symplectic and quadratic modules}\label{ss:modules}
Let $A$ be a finite abelian group.

\begin{defe}
A pairing $\omega\colon A\otimes A\ra \bCt$ is called \emph{skew} if $\omega(a,a)=0$ for all $a\in A$, and \emph{symmetric} if $\omega(a,b)=\omega(b,a)$ for all $a,b\in A$.
\end{defe}

Assume from now on that $\omega$ is either skew or symmetric.
If $L\subset A$ is a subgroup, write $L^\perp$ for the set
$$L^\perp:=\{a\in A\mid \omega(x,a)=1\,\forall x\in L\}.$$

\begin{defe}\label{def:isotropic}  A subgroup $L\subset A$ is \emph{isotropic} if $L\subset L^\perp$, \emph{coisotropic} if $L^\perp\subset L$, and \emph{Lagrangian} if $L=L^\perp$.
\end{defe}

\begin{defe}
The  \emph{kernel}  $\Ker(\omega)\subset A$ is defined by $\Ker(\omega):=A^\perp$.  A  skew or symmetric pairing is called  \emph{non-degenerate} if $\Ker(\omega)=0$.
\end{defe}

\begin{defe}
A  skew pairing is called \emph{sym\-plec\-tic} if it is non-degen\-erate.
We call a pair $(A,\omega)$ consisting of a finite abelian $p$-group $A$ and a symplectic (resp. non-degenerate symmetric) pairing $\omega$ a \emph{symplectic module} (resp. \emph{quadratic module}).
\end{defe}

Any skew $\omega$ induces a symplectic pairing on $A/\Ker(\omega)$ (and analogously if $\omega$ is symmetric). If $\omega$ is non-degenerate and $L\subset A$ is Lagrangian, then $|A|=|L|^2$.

 \subsection{The discriminant of a quadratic module} \label{sss:discriminant}
Let $(A,q)$ be a quadratic module. According to Corollary \ref{c:pairingDet}, $q$ defines an isomorphism $\det(q)\colon \det(A)\otimes \det(A)\isom \Fp$. If we choose a trivialization $\det(A)\isom \Fp$, then $\det(q)$ amounts to an isomorphism $\Fp\to\Fp$, i.e. to an element $\delta(q)\in\Fpt$. It is easy to show that the class of $\delta(q)\in \Fptt$ is independent of the chosen trivialization.
\begin{defe}We call $\delta(q)\in\Fptt$ \emph{the discriminant} of $(A,q)$.
\end{defe}

\begin{rem} When $A$ is a vector space over $\Fp$, this discriminant is often called the `signed discriminant' in the literature (e.g. in \cite{Lam}). It differs from the traditional discriminant of a quadratic form by a factor of $(-1)^{\frac{d(d-1)}{2}}$, where $d=\len A$, because of our convention  \eqref{dualityconvention}.
\end{rem}

\begin{exam} \label{e:angles}  Fix a non-trivial character $\psi_0\colon \Fp\to \bCt$. (It should agree with the one in \S\ref{sss:dualityAb}.) For $a\in \Fpt$, let
\[
\langle a \rangle\colon  \Fp\otimes \Fp\ra \bCt
\]
be the symmetric form $(x,y)\mapsto \psi_0(axy)$.\footnote{Note that the isomorphism class of $\langle a \rangle$ depends only on the image of $a$ in $\Fptt$.}
Set
\[
\langle a_1,\ldots,a_d \rangle:=\langle a_1 \rangle\oplus \ldots \oplus \langle a_d \rangle.
\]
Then \[
\delta(\langle a_1,\ldots,a_d \rangle)=(-1)^{\frac{d(d-1)}{2}}a_1a_2\ldots a_d.
\]
\end{exam}

\subsection{The Witt group} \label{ss:wittGroup}

\begin{ntn}\label{n:circ}
Given any quadratic module $(A,q)$, write $A^\circ$ for the same abelian group $A$ equipped with the pairing $q^\circ(a,b)=q(a,b)\inv$.
\end{ntn}

\begin{defe} \label{d:lagrangianCor}
 A \emph{Lagrangian correspondence}  $L\colon (A_1,q_1)\to (A_2,q_2)$ is a Lagrangian subgroup of $A_1^\circ\oplus A_2$.
 Quadratic modules $(A_1,q_1)$ and $(A_2,q_2)$ are \emph{Witt-equi\-va\-lent} if there exists a Lagrangian correspondence between them.
\end{defe}

\begin{exam} \label{ex:WittEquivalentGroups}
 If $\phi\colon A \to A'$ is an isomorphism of quadratic modules,  then the graph of $\phi$ in $A^\circ\oplus A'$ is a Lagrangian correspondence. Thus isomorphic quadratic modules are Witt-equivalent.
\end{exam}

\begin{prop} \label{p:WittEquivalenceRelation} Witt-equivalence is an equivalence relation.
\end{prop}

\begin{rem} In fact, the proof of Proposition \ref{p:WittEquivalenceRelation} in \S\ref{pf:equiv} defines a category whose objects are quadratic modules and whose morphisms are Lagrangian correspondences.
\end{rem}

\begin{defe} \label{d:WittGroup}
The Witt group $\bW$ is the group of Witt-equivalence classes of quadratic modules, where addition is given by the direct sum.
\end{defe}

\subsubsection{Description of the Witt group}  \label{sss:discriptionOfWittGroup}

Let $\bW_0$ be the group whose underlying set is
$\bW_0=\bZ/2\bZ \times \Fptt$ and whose multiplication is given by
\begin{equation}\label{eq:Qmult}
(e,d)\cdot(e',d'):=(e+e', (-1)^{ee'} dd').
\end{equation}

\begin{prop} \label{p:WittIsoQuadraticClasses}
The map $(\len, \delta)\colon  \bW\ra \bW_0$ is an isomorphism of groups.\footnote{It is easy to see that $\len\colon \bW\ra \bZ/2\bZ$ and $\delta\colon  \bW\ra \Fptt$ are well-defined functions.} Its inverse $\qm$ is given by
\[
\qm(1,a)= \langle a \rangle \qquad  \qm(0,a)= \langle -1,a \rangle.
\]
\end{prop}
The proof in \S\ref{pf:WittIso} reduces this proposition to the case of quadratic forms over finite fields studied, for example, in \cite[\S II.2]{Lam}.

\begin{rem} \label{r:WittEquivalenceWithDeterminant}
In particular, the quadratic module $q$ is Witt-equivalent to
$\langle \delta(q) \rangle$ if $\len(q)$ is odd, and to  $\langle -1, \delta(q) \rangle$ if $\len(q)$ is  even.
\end{rem}

\subsection{The Maslov index}\label{ss:Maslov}
Let $(A,\omega)$ be a symplectic module.  Let $L=(L_1,\ldots,L_m)$ be a sequence of Lagrangians in $A$. For notational convenience, set $L_{m+1}:=L_1$. Consider the complex of abelian groups
\[C_L:=
\xymatrix@C0.5in{[\bigoplus_{i=1}^m L_i\cap L_{i+1} \ar[r]^{\partial} &  \bigoplus_{i=1}^m L_i \ar[r]^{\sum} & A]}
\]
where $\Sigma$ is the summation, and $\partial$ maps $L_i\cap L_{i+1}$ into $L_i\oplus L_{i+1}$ by $a\mapsto (a,-a)$.
Let $T_L$ be the homology of $C_L$ at the center term.
Define a pairing $q_L\colon T_L\ten T_L\to\bCt$ by
\begin{equation}\label{eq:MDef} q_L(a,b):=\prod_{m\geq i>j\geq 1} \omega(a_i,b_j).\end{equation}

\begin{prop}[\cite{TeruMaslov}]
\label{p:MaslovQuadraticForm}
The pairing $q_L$ is well-defined, symmetric, and non-degenerate.
\end{prop}

\begin{defe} \label{d:MaslovIndex}
The {\it Maslov index} of $L$ is the class of the quadratic module $(T_L,q_L)$ in $\bW$.
\end{defe}

\begin{prop}[cf. \cite{TeruWeil}, \S 4.2]
\label{p:dimOfMaslov}\hfill
\begin{enumerate}
\item The length of $T_L$ equals
\[
\frac{m-2}{2} \len(A) - \sum_{1\leq i \leq m} \len(L_i\cap L_{i+1}) +2 \len(\cap_{i=1}^m L_i).
\]
\item Equip the $L_i$ with arbitrary orientations and let $\tL_i$ denote the resulting oriented Lagrangians (see Definition \ref{d:oriented}). Then
\[
\delta(q_L)=(-1)^{\frac{1}{2} \sum_{i\neq j} m_im_j} \prod_{i=1}^m \theta(\tL_i,\tL_{i+1})
\]
where $m_i:=\len(L_i/L_i\cap L_{i+1})$.
\end{enumerate}
\end{prop}

\begin{thm}[\cite{TeruMaslov}] \label{t:MaslovIsCocycle}
If $L_1,\ldots,L_m$ are Lagrangians in $(A,\omega)$, then
\[\tau(L_1,\ldots,L_m)=\tau(L_2,\ldots,L_m,L_1)=-\tau(L_m,L_{m-1},\ldots,L_1)\]
and
\[\tau(L_1,\ldots,L_m) = \tau(L_1,\ldots,L_k) + \tau(L_1,L_{k},\ldots,L_m)\]
for any $k<m$.
\end{thm}

\subsubsection{The Maslov cocycle is a coboundary on oriented Lagrangians} \label{sss:MaslovIsCoboundary}

Let $(A,\omega)$ be a symplectic module.
Let $\tL_1,\tL_2$ be oriented Lagrangians in $A$. Set
\begin{equation}\label{def:Theta}
\Theta(\tL_1,\tL_2)=\qm({\len(L_1/L_1\cap L_2)}, \theta(\tL_1,\tL_2))\in \bW,
\end{equation}
where  $\theta$ is the relative orientation (\S\ref{sss:relativeOrientLagrangians}), and $\qm$ was defined in Proposition \ref{p:WittIsoQuadraticClasses}. The following result is a restatement of ~\cite{Indians}, Prop. 2.1, in the finite group setting; it follows easily from  \eqref{eq:Qmult} and Proposition \ref{p:dimOfMaslov}.

\begin{prop} \label{p:MaslovIsCoboundary}
Let $\tL_i=(L_i, u_i)$, $1\leq i\leq m$ be oriented Lagrangians. Then we have
\[
\tau(L_1,\ldots,L_m)=\sum_{i\in \bZm} \Theta(\tL_i,\tL_{i+1}).
\]
\end{prop}

\subsection{The $\gamma$-index (after A. Weil)} \label{ss:GammaIndex}

\begin{defe} \label{d:GammaIndex}
For any quadratic module $(A,q)$, define
\begin{equation}\label{e:defgamma}\gamma(q)=\frac{1}{\sqrt{|A|}}\sum_{a\in A} q(a,a/2).\end{equation}
For any $a\in\Fpt$, write $\gamma(a)$ for $\gamma(\langle a\rangle)$ (in the notation of Example \ref{e:angles}). \end{defe}

\begin{rem} In \eqref{e:protogamma} we chose  $\psi_0(x)=\exp(2\pi i x/p)$ for concreteness.
\end{rem}

\begin{prop}\label{p:gammaIsHomomorphism}
This $\gamma$ defines a homomorphism  $\gamma\colon  \bW\to\bCt$.
\end{prop}

The proof is given in \S\ref{pf:gammaIsHomomorphism}. Since $|\bW|=4$ by Proposition \ref{p:WittIsoQuadraticClasses}, we obtain:

\begin{cor} \label{c:GammaOfMinus}
For any quadratic module $(A,q)$, $\gamma(q)$ is a fourth root of unity.
\end{cor}

Combining Remark \ref{r:WittEquivalenceWithDeterminant} and Proposition \ref{p:gammaIsHomomorphism} we obtain:

\begin{cor} \label{c:GammaAndDiscriminant}
Let $(A,q)$ be a quadratic module. Then\footnote{Note that $(d-1)^2\equiv 0\bmod 4$ if $d$ is odd, and $(d-1)^2\equiv 1\bmod 4$ if $d$ is even.}
\[
\gamma(q)= \gamma(1)^{-(d-1)^2} \gamma(\delta(q)).
\]
\end{cor}

Finally, the equality $\langle a,b\rangle\>=\langle 1,ab\rangle$ in $\bW$ gives the useful result
\begin{cor} \label{c:GammaOfProducts}
For any  $a,b\in \Fpt$,
$
\gamma(a)\gamma(b)=\gamma(1)\gamma(ab).
$
\end{cor}

\begin{rem}\label{r:GammaAsGaussSum}\emph{(Explicit values)} Recall that the definition of $\langle a \rangle$ in Example \ref{e:angles} depends on the choice of an additive character $\psi_0$.  For some $m$ coprime to $p$ we have  $\psi_0(x)=e^{2\pi im x/p}$, and
one can show that $\gamma(a)=\left( \frac{am}{p} \right)\epsilon_p$, where
$\left(\frac{am}{p}\right)$ is the Legendre symbol and
\[\epsilon_p=\begin{cases}
 -1 & p\equiv 1 \mod 4 \\
 -i & p\equiv 3 \mod 4.
\end{cases}
\]
\end{rem}

\section{Proofs}\label{s:proofs}

\subsection{Lemma \ref{l:chainOfPolarizations}}\label{pf:chainOfPolarizations}

Let us call a chain of consecutively neighboring polarizations, as in the statement of the Lemma, simply `a chain.'

We proceed by induction on the cardinality of $\fg$. Set $\fg_0:=Z(\fg)\cap\ker(f)$. Then $\fg_0$ is an ideal of $\fg$ contained in every polarization of $f$. Moreover, $f$ descends to a character of the Lie ring $\fg/\fg_0$, for which $\fp_1/\fg_0$ and $\fp_2/\fg_0$ are polarizations. If $\fg_0\neq 0$, then by induction we obtain a chain in $\fg/\fg_0$, and it lifts to a chain in $\fg$. We can now, therefore, assume that $\fg_0=0$.

In particular, $\fg$ is  non-abelian (if not 0), so we can choose $x\in\fg$ such that $0\neq[\fg,x]\subset Z(\fg)$. Using the skew form $(a,b)\mapsto f([a,b])$ on $\fg$ (see \S\ref{ss:modules}), let $\fh=\{x\}^\perp$ . Then $\fh$ is a coisotropic ideal in $\fg$.

Define $\fp'_i=\fp_i\cap\fh+\fh^\perp$. It is easy to check that  $\fp'_i$ is a polarization and  that it is a neighbor of $\fp_i$. It follows that $\fp'_1$ and $\fp'_2$ are polarizations for the restriction of $f$ to $\fh$, and, since $\fh$ is strictly smaller than $\fg$, they are joined by a chain in $\fh$, so \textit{a fortiori} in $\fg$. In total, we have constructed a chain $\fp_1,\fp'_1,\ldots,\fp'_2,\fp_2$ as desired.

\subsection{Lemma \ref{p:symfla}}\label{pf:symfla}
 Let  $\delta\in V_{f,\fp_1}$ be supported (as a function on $G$) on $P_1$ and equal there to $\delta(g)=\chi_f(g)$. We have
\begin{equation*}\begin{aligned}
\alpha(\fp_1,\fp_2,\fp_3)&=
\Phi_{\fp_1,\fp_2}\circ \Phi_{\fp_2,\fp_3}\circ \Phi_{\fp_3,\fp_1}\delta(1)\\
&=N_0(\fp_1,\fp_2,\fp_3)\sum_{p_3,p_2,p_1}
(\chi_f(p_3) \chi_f(p_2) \chi_f(p_1))^{-1}\delta(p_3p_2p_1)\end{aligned}\end{equation*}
where $N_0(\fp_1,\fp_2,\fp_3):= {\left(|\fp_1||\fp_2||\fp_3||\fp_1\cap\fp_2||\fp_2\cap\fp_3||\fp_3\cap\fp_1|\right)^{-1/2}}.$
The sums are over all $p_i\in P_i$, but the summand is only non-zero when
$p'_1:=(p_3p_2)^{-1}$ lies in $P_1$. Writing $\delta(p_3p_2p_1)=\chi_f(p')^{-1}\chi_f(p_1)$ we
obtain simply
\begin{equation*}\begin{aligned}
\alpha(\fp_1,\fp_2,\fp_3)&=N_0(\fp_1,\fp_2,\fp_3)\sum_{p_1\in P_1, (p_3,p_2,p'_1)\in S}
(\chi_f(p_3)\chi_f(p_2) \chi_f(p'_1))^{-1}.\end{aligned}\end{equation*}
The summand is independent of $p_1$, thus introducing a factor of $|\fp_1|=\sqrt{|\fg||\Ker\,B_f|}$ (cf. Remark \ref{r:sizeOfPolarization}) to yield
the desired formula \eqref{eq:form}.

\subsection{Lemma \ref{l:neighborfla}}\label{pf:neighborfla}
In \eqref{eq:form} we have $p_3=(p_2p_1)^{-1}$. Looking at the Campbell-Hausdorff series, the fact that $\fp_1$ and $\fp_2$ are neighboring polarizations implies that
\begin{equation}\label{qqq}-\log p_3\equiv\log p_1 + \log p_2 + \tfrac12[\log p_2,\log p_1] \mod \fp_1\cap\fp_2\cap\ker f.\end{equation}
In particular, $-\log p_3-\log p_1$ lies in $\fp_2$, and
$$(p_3,p_2,p_1)\mapsto(q_3,q_2,q_1):=(\log p_3,-\log p_3-\log p_1,\log p_1)$$
is a well-defined bijection $S\to\fS$.
We use this bijection to rewrite \eqref{eq:form} as a sum over $\fS$. Using \eqref{qqq} we calculate that  the summand
\[(\chi_f(p_3) \chi_f(p_2)\chi_f(p_1))^{-1}=f(\tfrac12[\log p_2,\log p_1]) =
f(\tfrac12[q_2,q_1])\]
so we obtain
\begin{equation}\label{eq:form2}\alpha(\fp_1,\fp_2,\fp_3)={N(\fp_1,\fp_2,\fp_3)}\sum_{(q_3,q_2,q_1)\in \fS} f(\frac 12[q_2,q_1]).\end{equation}
Now it is easy to verify that the summand depends only on the class of $(q_3,q_2,q_1)$ in $\fS/\fR$. That class has size
$$\frac{|\fp_1\cap\fp_2||\fp_2\cap\fp_3||\fp_3\cap\fp_1|}{ |\fp_1\cap\fp_2\cap\fp_3|}=\frac{N'(\fp_1,\fp_2,\fp_3)}{N(\fp_1,\fp_2,\fp_2)},$$ so we obtain \eqref{eq:cocycleNeighbors}.

\subsection{Proposition \ref{p:cocycleIsMaslov}}\label{pf:cocycleIsMaslov}
Since all polarizations in a Heisenberg Lie ring are neighbors (\S\ref{sss:neighboringLagrangians}), we may use Lemma \ref{l:neighborfla}. The order-reversing map
$$F\colon L_1\oplus L_2\oplus L_3 \to \fp_3\oplus \fp_2\oplus \fp_1$$
defines a linear isomorphism from $T_L$ (\S\ref{ss:Maslov}, with $L:=(L_1,L_2,L_3)$) to $\fS/\fR$. Using this to rewrite \eqref{eq:cocycleNeighbors} as a sum over $T_L$, we find
$$\alpha(\fp_1,\fp_2,\fp_3)=N'(\fp_1,\fp_2,\fp_3) \sum_{a\in T_L}
 q_L(a,a/2).$$
 It only remains to compare this with the definition \eqref{e:defgamma} of $\gamma$.

\subsection{Proposition \ref{p:intertwinersforHeisenberg}}\label{pf:intertwinersforHeisenberg}
By Propositions \ref{p:MaslovIsCoboundary} and \ref{p:gammaIsHomomorphism}, we have
$$
\gamma(\tau(L_1,L_2,L_3))=\gamma(\Theta(\tL_1,\tL_2))\gamma(\Theta(\tL_2,\tL_3))
\gamma(\Theta(\tL_3,\tL_1)).
$$
Therefore, it remains to check that
$$\gamma(\Theta(\tL_1,\tL_2))=\beta(\tL_1,\tL_2):=\gamma(1)^{-m(L_1,L_2)^2} \gamma(\theta(\tL_1,\tL_2)).$$
But this follows immediately from Corollary \ref{c:GammaAndDiscriminant} and the definition \eqref{def:Theta} of $\Theta$.

\subsection{Proposition \ref{p:WittEquivalenceRelation}}\label{pf:equiv}
Suppose given Lagrangian correspondences
$L_{21}\subset A_1^\circ\oplus A_2$ and
$L_{32}\subset A_2^\circ\oplus A_3$. Define
$L_{31}\subset A_1^\circ\oplus A_3$ by
\[
L_{31}:=\{(a,c)\mid\exists b\in A_2\mbox{ with } (a,b)\in L_{21}\mbox{ and } (b,c)\in L_{31}\}.
\]
Then $L_{31}$ is Lagrangian.

\subsection{Proposition \ref{p:WittIsoQuadraticClasses}}\label{pf:WittIso}
The map $(\len,\delta)$ is surjective, since $\qm$ is a right inverse. So it will suffice to show that $\bW$ has four elements.

Suppose $(A,q)$ is a quadratic module, and $n$ is the smallest integer such that $p^nA=0$. Then $A_{n-1}:=p^{n-1}A$ is an isotropic subgroup of $A$, and $A$ is Witt-equivalent to $A_{n-1}^\perp/A_{n-1}$:
$$\{(a,b)\in A^\circ\oplus  A_{n-1}^\perp/A_{n-1}\mid a\equiv b \bmod{A_{n-1}}\}$$
is a Lagrangian correspondence. Iterating this process, we see that every quadratic module is  Witt-equivalent to a $p$-torsion one,  i.e. to a vector space over $\Fp$.
A symmetric form on a vector space can be diagonalized; this shows that $\bW$ is generated by $\langle 1\rangle$ and $\langle a \rangle$ where $a\in \Fpt$ is not a square.

If $-1$ is a square in $\Fpt$, then $\langle 1\rangle=\langle-1\rangle=-\langle1\rangle$ in $\bW$, so that $\langle1, 1\rangle=0$ (and, similarly, $\langle a,a\rangle=0$) in $\bW$;  thus $\bW$ is isomorphic to $\bZ/2\bZ\times \bZ/2\bZ$.

If $-1$ is not a square, then $\langle a\rangle=\langle -1\rangle=-\langle 1\rangle$ in $\bW$, and it remains to show that $\langle1,1,1,1\rangle=0$. Every element of a finite field is a sum of two squares; write  $-2=x^2+y^2$. Then the linear span of $(x,y,1,1)$ and $(1,1,-x,-y)$ in $\langle1,1,1,1\rangle$ is Lagrangian.

%
%
\subsection{Proposition \ref{p:gammaIsHomomorphism}}\label{pf:gammaIsHomomorphism} Suppose $(A_1,q_1),(A_2,q_2)$ are quadratic modules. One sees immediately that
\begin{equation}\label{additive}\gamma(q_1\oplus q_2)=\gamma(q_1)\gamma(q_2).\end{equation}
We must check that $\gamma(q_1)=\gamma(q_2)$ whenever $(A_1,q_1)$ and $(A_2,q_2)$ are Witt equivalent.
It suffices to show that $\gamma(q)=1$ whenever  $(A,q)$ is a Witt-equivalent to $0$.  For then equation \eqref{additive} gives
$$\gamma(q_1^\circ)\gamma(q_2)=\gamma(q_1^\circ\oplus q_2)=1.$$
In particular,
$\gamma(q_1^\circ)\gamma(q_1)=1$. Therefore $\gamma(q_1)=\gamma(q_2)$. It also follows from \eqref{additive} that $\gamma$ is a homomorphism.

Suppose, then that $(A,q)$ contains a Lagrangian $L$, and let $S\subset A$ be a set of representatives of $A/L$. We find
\[\begin{aligned}\gamma(q)&=\frac{1}{\sqrt{|A|}}\sum_{a\in L}\sum_{s\in S} q(a+s,(a+s)/2) \\
& = \frac{1}{\sqrt{|A|}} \sum_{a\in L}\sum_{s\in S} q(a,a/2)q(a,s)q(s,s/2) \\
& = \frac{1}{\sqrt{|A|}}  \sum_{s\in S}q(s,s/2)\sum_{a\in L}q(a,s) \\
& = q(0,0)=1.\end{aligned}
\]
The first equality is the definition, the second is bi-additivity, the third is because $L$ is isotropic, and the fourth uses that
$\sum_{a\in L}\psi(q(a,s))=0$
unless $s\in L$.

\bibliography{intertwiners}

\begin{thebibliography}{Kam05}

\bibitem[BD06]{DB}
D.~Boyarchenko and V.~Drinfeld.
\newblock A motivated introduction to character sheaves and the orbit method
  for unipotent groups in positive characteristic.
\newblock {\em arXiv:math/0609769v1}, 09/27/2006.

\bibitem[Del87]{DeligneDet}
P.~Deligne.
\newblock Le determinant de la cohomologie.
\newblock {\em Contem. Math.}, 67:93--177, 1987.

\bibitem[G\'77]{Gerardin}
P.~G\'{e}rardin.
\newblock Weil representation associated to finite fields.
\newblock {\em J. Algebra}, 46(1):54--101, 1977.

\bibitem[GH08a]{RonnyGeometricWeil}
S.~Gorevich and R.~Hadani.
\newblock Geometric {W}eil representation.
\newblock {\em arXiv:math/0610818v2}, 2008.

\bibitem[GH08b]{RonnyCanonical}
S.~Gorevich and R.~Hadani.
\newblock Quantization of symplectic vector spaces over finite fields.
\newblock {\em arXiv:0705.4556v3}, 2008.

\bibitem[How73]{Howe-Character}
Roger~E. Howe.
\newblock On the character of {W}eil's representation.
\newblock {\em Trans. Amer. Math. Soc.}, 177:287--298, 1973.

\bibitem[How77]{Howe-Nilpotent}
Roger~E. Howe.
\newblock On representations of discrete, finitely generated, torsion-free,
  nilpotent groups.
\newblock {\em Pacific J. Math.}, 73(2):281--305, 1977.

\bibitem[Kam05]{MasoudWeil}
M.~Kamgarpour.
\newblock Weil representations over finite fields.
\newblock {\em www.masoudkamgarpour.com/Media Files/masterthesis.pdf}, 2005.

\bibitem[Kir62]{Kirillov}
A.~A. Kirillov.
\newblock Unitary representations of nilpotent lie groups.
\newblock {\em Uspehi Mat. Nauk}, 17:57--110, 1962.

\bibitem[Lam05]{Lam}
T.Y. Lam.
\newblock {\em Introduction to quadratic forms over fields}, volume~67 of {\em
  Graduate Studies in Mathematics}.
\newblock Amer. Math. Soc., 2005.

\bibitem[Lan05]{LA}
S.~Lang.
\newblock {\em Algebra}, volume 211 of {\em Graduate Text in Mathematics}.
\newblock Springer-Verlag, 2005.

\bibitem[Laz54]{Lazard}
M.~Lazard.
\newblock Sur les groupes nilpotents et anneaux de {L}ie.
\newblock {\em Ann. Sci. Eco. Norm. Sup.}, (31) 71:101--190, 1954.

\bibitem[LP80]{LP}
G.~Lion and P.~Perrin.
\newblock Extensions des repr\'{e}sentation de groupes unipotents p-adic.
\newblock {\em Lect. Notes in Math.}, 880, 1980.

\bibitem[PPS00]{Indians}
R.~Parimala, R.~Preeti, and R.~Sridharan.
\newblock Maslov index and a central extension of the symplectic group.
\newblock {\em K-Theory}, 19:29--45, 2000.

\bibitem[Qui73]{Quillen}
D.~Quillen.
\newblock Higher algebraic {K}-theory {I}.
\newblock {\em Lecture notes in Mathematics}, 341:85--143, 1973.

\bibitem[RY02]{Reichstein2002}
Z.~Reichstein and B.~Youssin.
\newblock A birational invariant for algebraic group actions.
\newblock {\em Pacific J. Math.}, 204(1):223--246, 2002.

\bibitem[Tho06]{TeruMaslov}
T.~Thomas.
\newblock The {M}aslov index as a quadratic space.
\newblock {\em Math. Res. Lett.}, 5-6:985--999, 2006.

\bibitem[Tho07]{TeruWeil}
T.~Thomas.
\newblock The character of the {W}eil representtion.
\newblock {\em arXiv: math/0610644}, 2007.

\bibitem[Wei64]{Weil}
A.~Weil.
\newblock Sur certains groups d'operateurs unitaires.
\newblock {\em Acta Math.}, 111:143--211, 1964.

\end{thebibliography}

\end{document}